\pdfoutput=1
\RequirePackage{ifpdf}
\ifpdf 
\documentclass[pdftex]{sigma}
\else
\documentclass{sigma}
\fi

\numberwithin{equation}{section}

\newtheorem{Theorem}{Theorem}[section]
\newtheorem*{Theorem*}{Theorem}
\newtheorem{Corollary}[Theorem]{Corollary}
\newtheorem{Lemma}[Theorem]{Lemma}

 { \theoremstyle{definition}

\newtheorem{Example}[Theorem]{Example}
\newtheorem{Remark}[Theorem]{Remark} }

\newcommand{\C}{{\mathbb C}}
\DeclareMathOperator{\End}{End}
\DeclareMathOperator{\dimension}{dim}
\DeclareMathOperator{\tr}{tr}
\DeclareMathOperator{\id}{id}
\DeclareMathOperator{\inv}{inv}
\newcommand{\Uqgl}{{U_q(\mathfrak{gl}(N+1))}}

\newcommand{\zz}{{\mu}}
\newcommand{\ii}{{\mathbf{i}}}
\newcommand{\jj}{{\mathbf{j}}}

\newcommand{\zo}{{(\mu)}}

\begin{document}

\allowdisplaybreaks

\newcommand{\arXivNumber}{1805.04884}

\renewcommand{\PaperNumber}{036}

\FirstPageHeading

\ShortArticleName{Explicit Central Elements of $U_q(\mathfrak{gl}(N+1))$}

\ArticleName{Explicit Central Elements of $\boldsymbol{U_q(\mathfrak{gl}(N+1))}$}

\Author{Jeffrey KUAN~$^{\rm a}$ and Keke ZHANG~$^{\rm b}$}

\AuthorNameForHeading{J.~Kuan and K.~Zhang}

\Address{$^{\rm a)}$~Department of Mathematics, Texas A\&M University, Mailstop 3368, College Station,\\
\hphantom{$^{\rm a)}$}~TX 77843-3368, USA}
\EmailD{\href{mailto:jkuan@tamu.edu}{jkuan@tamu.edu}}

\Address{$^{\rm b)}$~Perimeter Institute, 31 Caroline St. N, Waterloo, ON, N2L 2Y5, Canada}
\EmailD{\href{mailto:kzhang@perimeterinstitute.ca}{kzhang@perimeterinstitute.ca}}

\ArticleDates{Received August 25, 2022, in final form May 19, 2023; Published online June 03, 2023}

\Abstract{By using Drinfeld's central element construction and fusion of $R$-matrices, we construct central elements of the quantum group $U_q(\mathfrak{gl}(N+1))$. These elements are explicitly written in terms of the generators.}


\Classification{16T05} 	
	
\section{Introduction}	

Several recent works (e.g., \cite{CGRS2,CGRS,KuanQRWQG,KuanJPhysA,KIMRN}), motivated by mathematical physics and probability, have used explicit central elements of various quantum groups. For the applications in those paper, the central elements need to be explicitly written in terms of the generators of the quantum group. In particular, in a probabilistic perspective, the probabilities are required to be non-negative, which requires explicit expressions for the central elements for calculations. Furthermore, in an asymptotic analysis, the analysis tools also require explicit expressions for the central elements.

The explicit generators of the center of $\Uqgl$ was first constructed for generic $q$
in~\cite{FaddeevRT} without proof (see~\cite{Arnaudon:1993yg}). Previous work of~\cite{ZGBCMP} applies Drinfeld's central element construction~\cite{Drin90} to universal $R$-matrices in order to construct central elements of quantum groups and to determine their eigenvalues on irreducible highest weight modules.

Jimbo in \cite{Jimbo1986} provided a method to construct a polynomial $z(x) = \sum_{i=0}^{N+1} z_ix^i$ in variable~$x$ such that all $z_i$’s are central. Tanisaki in \cite{Tanisaki1992KILLINGFH} proved that all these $z_i$’s generate the whole centre by quantized Harish-Chandra isomorphism. However, the expression of $z_i$’s in \cite{Jimbo1986} are not explicit enough for applications cited in the first paragraph.

By using Jimbo's \cite{Jimbo1986} formula for the $R$-matrix of $\Uqgl$, further work by the same authors~\cite{GZB} explicitly writes a quantum Casimir element of $\Uqgl$ with a formula for its eigenvalues. Using an explicit formula for the universal $R$-matrix in~\cite{KirillovResh}, the authors~\cite{ZGBJPhysA} write an explicit (but somewhat complicated) expression for general Casimir elements in quantum groups.

In this paper, we apply Drinfeld's central element construction to the fused $R$-matrices of $\Uqgl$ in~\cite{Jimbo1986}, rather than the universal $R$-matrices of~\cite{KirillovResh}. The resulting central elements appear to be slightly simpler than the previous expressions. The proof requires some new ingredients, notably relations between the root vectors~\cite{KT91,Lusztig90,Xi1994} and some elementary knowledge of coset representatives of symmetric groups.

We also note the work in~\cite{Junbo}, which shows that the Casimir elements (and a trivial central element) generate the entire center of $\Uqgl$, although this had been known as early as~\cite{Rosso87}. Additionally, the paper~\cite{LXZ} explicitly writes two algebraically independent central elements in~$\Uqgl$. Recently, Dai~\cite{Dai} provided the explicit generators of the center of all (simply-connected) quantum groups of finite type by following~\cite{ZGBCMP}. In principle, it should be possible to match those central elements to the ones here, but we do not pursue this direction here.

\section{Notations and backgrounds}
\subsection{Symmetric groups}
Recall that the symmetric group is the Weyl group of the special linear group. Define the usual action of the symmetric group $S_m$ on $\mathbb{N}^m$ by $ \sigma(x_1,\dots,x_m)= (x_{\sigma(1)}, \dots, x_{\sigma(m)})$. For any $A\in \mathbb{N}^m$, define $ H_A \le S_m $ to be the subgroup $ \{\sigma \in S_m\colon \sigma(A)=A\}$. Let $ D_A \subset S_m $ be the set of (left) coset representatives of $H_A$ with the fewest inversions: in other words, $ \sigma \in D_A $ if and only if $ \inv (\sigma) \le \inv (\tau) $ for every $ \tau \in \sigma H_A$. Here, $\inv$ is the length function defined as the fewest number
of simple reflections constituting the Weyl group element.

 \begin{Example} Take $ A= (3,3,3,2,2,1).$ Then $H_A=S_3 \times S_2 \times S_1 \leq S_6$. We have $ (34) \in D_A $ with $ (34) \cdot A= (3,3,2,3,2,1)$, and $ (35) \in D_A $ with $ (35)\cdot A=(3,3,2,2,3,1) $. However, $ (134)\cdot A= (34) \cdot A$ and $\inv((134)) > \inv((34))$, so $ (134) \notin D_A $.
 \end{Example}

 We recall (see, e.g., \cite{Carter}) that each coset of $H_A$ has a unique representative $\sigma \in D_A$, and that $\inv(\sigma\tau) = \inv(\sigma) + \inv(\tau)$ for every $\tau \in H_A$.

 \begin{Lemma}\label{Gene}
 	Suppose that $\tau,\tau' \in D_A$. Then there exist a sequence of elements $\tau_0,\dots,\tau_l\in D_A$ such that $\tau_0 = \tau'$, $\tau_l=\tau$ and every $\tau_{j+1}\tau_j^{-1}$ is a transposition for $j=0,\dots,l-1$.
 \end{Lemma}
 \begin{proof}
 	It suffices to prove this statement when either $\tau$ or $\tau'$ is the identity permutation $e$, because the two sequences can be concatenated. So suppose that $\tau$ is an arbitrary element of $D_A$ and $\tau'=e$. Let $s_k\dots s_1$ be a minimal word representation of $\tau$ and set $\tau_j = s_j \dots s_1$. Then $\tau_{j+1}\tau_j^{-1}=s_{j+1}$, which is a transposition. So it remains to show that $\tau_j \in D_A$. If it were not, then there would exist a transposition $s\in H_A$ such that $\inv(\tau_js) = \inv(\tau_j)-1$. But then $\inv(\tau s) = \inv(\tau)-1$, contradicting the assumption that $\tau \in D_A$.
 	
 	Now suppose that $\tau=e$ and $\tau'$ is an arbitrary element of $D_A$. By the previous paragraph, there exist a sequence of elements $\tilde{\tau}_0=e,\tilde{\tau}_1,\dots,\tilde{\tau}_l=\tau'$ in $D_A$ such that every $\tilde{\tau}_{j+1}\tilde{\tau}_j^{-1}$ is a~transposition. Setting $\tau_j = \tilde{\tau}_{l-j}$, we have that $\tau_0=\tau',\dots,\tau_l=e$ is a sequence of elements in $D_A$ and ${\tau}_{l-j-1}{\tau}_{l-j}^{-1}$ is a transposition for every $j$. The latter equality is equivalent to the condition that for every $k$, $\tau_k = s\tau_{k+1}$ for some transposition $s$. Since this is also equivalent to the condition that $\tau_{k+1}\tau_k^{-1}$ is a transposition for every $k$, this finishes the proof.
 \end{proof}

 Let $ \mathcal{B}_m^{(N)} $ denote the set of sequences $ \mu=(\mu_0,\dots,\mu_N) $ of non-negative integers such that $ \mu_0+\dots+\mu_N=m $. For any $ \mu \in \mathcal{B}_m^{(N)} $, let $ H^\mu \le S_m$ denote the subgroup $ S_{\mu_0}\times \dots \times S_{\mu_N} $, and likewise let $ D^\mu $ denote the set of left coset representatives with the fewest inversions. Define $\mathcal{B}_m$ to be the union $\bigcup_{N\geq 1} \mathcal{B}_m^{(N)}$.

Let $\mathcal{W}_m \subset \mathbb{N}^m$ denote the subset of elements $(i_1,\dots,i_m)$ satisfying $i_1 \leq \dots \leq i_m$. For $\ii \in \mathcal{W}_m$, and assuming that $N\geq i_m$, let $\mu^{(N)}(\ii) \in \mathcal{B}_m^{(N)}$ be defined by
\[
\mu^{(N)}(\ii)_k = \vert \{ l \in \{1,\dots,m\}\colon i_l=k \} \vert.
\]
For $\ii \in \mathcal{W}_m$ satisfying $i_m \leq N$, we have a natural isomorphism between the subgroups $H_{\ii}$ and $H^{\mu^{(N)}(\ii)}$. Thus there is also a natural bijection between $D_{\ii}$ and $D^{\mu^{(N)}(\ii)}$.

Given $\ii \in \mathcal{W}_m$, define the equivalence relation $\sim_{\ii}$ on $S_m$ by
\[
\tau \sim_{\ii} \sigma \text{ if and only if } \sigma^{-1}\tau \in H_{\ii}.
\]
In words, $\tau \sim_{\ii} \sigma$ if and only if they are in the same left coset of $H_{\ii}$. For any $\tau \in S_m$, there is a unique $\sigma \in D_{\ii}$ such that $\tau \sim_{\ii} \sigma$. Define $d_{\ii}(\tau)= \inv\big(\sigma^{-1}\tau\big)$. In other words, if $\tau$ is written uniquely as $\tau=\sigma \xi$ for $\sigma \in D_{\ii}$ and $\xi \in H_{\ii}$, then $d_{\ii}(\tau) = \inv(\xi)$. We will also let $\sigma_{\ii}(\tau)$ and $\xi_{\ii}(\tau)$ denote the two permutations in the unique expression $\tau = \sigma \xi$.

Finally, given $\tau \in S_m$, let $\bar{\tau}$ denote the reversed permutation $\bar{\tau}(k)= \tau(m+1-k)$.

We conclude this section by noting the following identity.
\begin{Lemma}\label{Leem}
For $\ii \in \mathcal{W}_m$, set $\mu = \mu^{(N)}(\ii)$. Then
\[
-N\mu_{0} -(N-2)\mu_{1} + \dots + N\mu_{N} = -Nm + 2(i_1+ \dots + i_m).
\]
\end{Lemma}
\begin{proof}
By definition,
\[
\mu_k = | l \in \{1,\dots,m\}\colon i_l=k |.
\]
We thus rewrite
\begin{align*}
-N\mu_{0} -(N-2)\mu_{1} + \dots + N\mu_{N} &= -N(\mu_0 + \dots + \mu_N) + 2\mu_1 + 4\mu_2 + \dots + 2N\mu_N \\
&= -Nm+ 2\mu_1 + 4\mu_2 + \dots + 2N\mu_N
\end{align*}
and note that
\[
2(i_1 + \dots + i_m) = 2\sum_{k=1}^N k\cdot | l \in \{1,\dots,m\}\colon i_l=k | = 2\sum_{k=1}^N k\cdot \mu_k,
\]
which shows the identity.
\end{proof}

\subsection{Quantum groups}
We use the following notation modified from \cite{Jimbo1986}.
Define $ U_q(\mathfrak{sl}(N+1)) $ to be the associative algebra over $ \mathbb{C} $ generated by the symbols $ q^{\pm h_i/2}$, $\hat{e}_{\pm, i}$, $1 \leq i\leq N$, under the following relations:
\begin{gather*}
q^{h_i/2} \cdot q^{-h_i/2}= q^{-h_i/2} \cdot q^{h_i/2}=1, \qquad q^{h_i/2} \cdot q^{h_{i'}/2}= q^{h_{i'}/2} \cdot q^{h_i/2},\nonumber\\
 q^{h_i/2} \hat{e}_{\pm, i'} q^{-h_i/2}= q^{\pm a_{ii'}/2} \hat{e}_{\pm, i'},\nonumber\\
[\hat{e}_{+,i}, \hat{e}_{-,i'}] =\delta_{i,i'} \frac{q^{h_i}-q^{-h_i}}{q-q^{-1}}, \nonumber\\
 \hat{e}_{\pm, i}\hat{e}_{\pm, i'} = \hat{e}_{\pm ,i'}\hat{e}_{\pm ,i}, \qquad |i-i'|\geq 2, \nonumber\\
 \hat{e}_{\pm, i}^2 \hat{e}_{\pm, i\pm 1} -\big(q+q^{-1}\big) \hat{e}_{\pm, i} \hat{e}_{\pm ,i\pm 1} \hat{e}_{\pm ,i} + \hat{e}_{\pm, i\pm 1} \hat{e}_{\pm, i}^2=0,\qquad 1\leq i, i\pm 1 \leq N. 
\end{gather*}
Here, $ (a_{i,i'})_{1\leq i,i' \leq N} $ denotes the Cartan matrix of type $ A_N $, i.e.,
\begin{equation*}
a_{ii'}=
\begin{cases}
\hphantom{-}2,& i = i',\\ -1, & i=i' \pm 1,\\ \hphantom{-}0,& \text{otherwise}.
\end{cases}
\end{equation*}
Then define $ U_q(\mathfrak{gl}(N+1)) $ by adjoining to $ U_q(\mathfrak{sl}(N+1)) $ the elements $ q^{\pm \epsilon_i/2}$, $0 \leq i \leq N $, so that $ q^{h_i}= q^{\epsilon_{i-1}-\epsilon_i} $ and that $ q^{\epsilon_0+\dots +\epsilon_N} $ belongs to the center.

The $m$-fold co-product is the algebra homomorphism \[\Delta^{(m)}\colon \ \Uqgl \rightarrow \underbrace{\Uqgl \otimes \dots \otimes \Uqgl}_{m} \]
such that \begin{gather}
 \Delta^{(m)}\big (q^{\pm \epsilon_i/2}\big)= q^{\pm \epsilon_i/2} \otimes \dots \otimes q^{\pm \epsilon_i/2},\qquad
 \Delta^{(m)} (\hat{e}_{\pm ,i} ) =\sum_{v=1}^{m} \hat{e}^{(v)}_{\pm ,i} ,\label{qerel}
\end{gather}
where
\[
\hat{e}^{(v)}_{\pm ,i} = \underbrace{q^{h_i/2} \otimes \dots \otimes q^{h_i/2}}_{v-1} \otimes \hat{e}_{\pm ,i} \otimes \underbrace{ q^{-h_i/2} \otimes \dots \otimes q^{-h_i/2}}_{m-v}.
\]
 We also have the reversed co-product
 \begin{gather*}
 \bar{\Delta}^{(m)} \big(q^{\pm \epsilon_i/2}\big)= q^{\pm \epsilon_i/2} \otimes \dots \otimes q^{\pm \epsilon_i/2}, \\
 \bar{\Delta}^{(m)} (\hat{e}_{\pm, i})= \sum_{v=1}^{m} \underbrace{q^{-h_i/2} \otimes \dots \otimes q^{-h_i/2}}_{v-1} \otimes \hat{e}_{\pm, i} \otimes q^{h_i/2} \otimes \dots \otimes q^{h_i/2}.
 \end{gather*}
We write $\Delta$, $\bar{\Delta}$ for $\Delta^{(2)}$, $\bar{\Delta}^{(2)}$. The map $ \Delta $ endows $ \Uqgl $ with a structure of a bi-algebra. It is also a Hopf algebra, but we will not need the counit and antipode.

Consider the following relations for $R$:
\begin{gather}
 R \Delta(u) = \bar{\Delta}(u) R \qquad \forall u\in \Uqgl. \label{Equi}
\end{gather}
 This admits a unique (up to a multiplicative constant) solution
$ R \in \End (V\otimes V) $, where $ V:=\mathbb{C}^{N+1} $ is the defining representation of $ \Uqgl $.

 Let us consider \eqref{Equi} in $ \Uqgl \otimes \End \big(\C^{N+1}\big) $. Then a family of solutions is given by (see \cite{Jimbo1986}) $R(x)= \sum_{0 \leq i,j \leq N} \hat{E}_{ij}'(x) \otimes e_{ji} $, where
 \begin{gather*} 
 \hat{E}'_{ij}(x) = \begin{cases}
 \big(x q^{(\epsilon_i+\epsilon_j-1)/2}\big)^{\mp 1} {E}'_{ij} , & i \lessgtr j,\\
 \big(xq^{\epsilon_i}- x^{-1}q^{- \epsilon_i}\big)/\big(q-q^{-1}\big), & i=j.
 \end{cases}
 \end{gather*}
 For this paper, only the value at $x=1$ will be used in the proofs, but for completeness Jimbo's result is stated in its full generality. Here, $ E'_{ij} $ are the root vectors defined recursively by
 \begin{gather*}
E'_{ij}= E_{ik}'E_{kj}' -q^{\pm 1} E_{kj}'E_{ik}',\qquad i \lessgtr k \lessgtr j, \qquad E'_{i-1,i} = \hat{e}_{+,i}, \qquad E'_{i,i-1} = \hat{e}_{-,i}
\end{gather*}
and $e_{ji}$ is the usual matrix which acts on the canonical basis $\{I_0,\dots,I_N\}$ of $\mathbb{C}^{N+1}$ by
\[
e_{ji}I_k = \delta_{ik} I_j.
\]

In \cite{GZB}, the authors define
\[
\hat{E}_{ij}
=
\begin{cases}
q^{(E_{ii}+E_{jj}-1)/2}E_{ij}, & i \neq j,\\
q^{E_{ii}}/\big(q-q^{-1}\big), & i=j,
\end{cases}
\]
where the modified root vectors are
\begin{equation*}
E_{ij}= E_{ik}E_{kj} -q^{- 1} E_{kj}E_{ik},\qquad i \lessgtr k \lessgtr j, \qquad E'_{i-1,i} = \hat{e}_{+,i}, \qquad E'_{i,i-1} = \hat{e}_{-,i}.
\end{equation*}
Using Jimbo's results, they show that
\begin{gather*}
R = \sum_{N \geq i\geq j \geq 0} \hat{E}_{ij} \otimes e_{ji} \in \Uqgl \otimes \End (V),\nonumber\\
R^T = \sum_{N \geq i\geq j \geq 0} \hat{E}_{ji} \otimes e_{ij} \in \Uqgl \otimes \End (V) 
\end{gather*}
satisfy
\begin{gather*}
R\Delta(u) = \bar{\Delta}(u) R,\qquad
R^T\bar{\Delta}(u) = {\Delta}(u) R^T.
\end{gather*}

We can solve \eqref{Equi} even more generally. Consider the $ m $-fold tensor of the defining representation $ V^{\otimes m}$, and let $P_mV^{\otimes m}$ be the symmetric projection. In other words, $P_mV^{\otimes m}$ is isomorphic to the module of highest weight~$m\nu$, where~$\nu$ is the highest weight of~$V$.

Then the solution of \eqref{Equi} in $ \Uqgl \otimes \End (P_m V^{\otimes m}) $ is given by the fused $R$-matrix~\cite{Jimbo1986}%
\begin{gather*} 
R^{}_{0,\{1,2,\dots,m\}}(x)= R_{0m}(x)R_{0m-1}(xq) \cdots R_{01}\big(xq^{ m-1}\big) P_m \in \Uqgl\otimes \End \big(P_m V^{\otimes m}\big) .
\end{gather*}

Therefore, the fused $R$-matrices
\begin{gather*}
R^{}_{0,\{1,2,\dots,m\}} = R_{0m}R_{0m-1}\cdots R_{01}P_m \in U_q(\mathfrak{gl}_{N+1})\otimes \End \big(P_m V^{\otimes m}\big), \\
R^{T}_{0,\{1,2,\dots,m\}} = R^T_{0m}R^T_{0m-1}\cdots R^T_{01}P_m \in U_q(\mathfrak{gl}_{N+1})\otimes \End \big(P_m V^{\otimes m}\big)
\end{gather*}
satisfy
\begin{gather}
R_{0,\{1,2,\dots,m\}} \Delta(u) = \bar{\Delta}(u) R_{0,\{1,2,\dots,m\} },\qquad
R^T_{0,\{1,2,\dots,m\}} \bar{\Delta}(u) = {\Delta}(u) R^T_{0,\{1,2,\dots,m\}}\label{RRel}
\end{gather}
in $\Uqgl \otimes \End (P_m V^{\otimes m})$ for all $u\in \Uqgl$.

Therefore,
\begin{equation*}
\Gamma_m: = R^T_{0,\{1,2,\dots,m\}} R_{0,\{1,2,\dots,m\}} \in \Uqgl \otimes \End \big(P_m V^{\otimes m}\big)
\end{equation*}
commutes with $ \Delta (u) $ for all $u\in \Uqgl$. By Drinfeld's central element construction \cite{Drin90}, the element
\begin{equation*}
C_m= \id \otimes\tr_q(\Gamma_m)
\end{equation*}
is central in $\Uqgl$, where the quantum trace $ \tr_q $ of an operator $A$ is defined by
\begin{equation*}
\tr_q(A) = \tr\big(q^{-2h_{\rho}}A\big) := \tr \big(q^{-N\epsilon_0-(N-2)\epsilon_1 - \dots +N\epsilon_N}A\big).
\end{equation*}

We also have the following relations between the root vectors:\footnote{Similar relations can be found in the paper \cite{Xi1994}, but there appear to be some typos. They can also be derived from \eqref{RRel} by applying $\id \otimes B$ to both sides, for suitable linear maps $B$ on $\mathrm{End}\big(P_2V^{\otimes 2}\big)$. One can also check that the relations hold in the explicit representations in Remark~\ref{Remark1} below.}
For $i<l$ and $j<k$,
\begin{gather}\label{comm}
{E}_{il}{E}_{jk} =
\begin{cases}
{E}_{jk}{E}_{il}, & i< j < k < l \text{ or } i<l<j<k,\\
q^{-1} {E}_{jk}{E}_{il}, & i=j<k<l, \\
q {E}_{jk}{E}_{il}, & i<j<k=l,\\
{E}_{jk}{E}_{il} + \big(q-q^{-1}\big) {E}_{jl}{E}_{ik}, & i<j<l<k,
\end{cases}
\\
\label{comm1}
{E}_{li}{E}_{kj} =
\begin{cases}
{E}_{kj}{E}_{li}, & i< j < k < l \text{ or } i<l<j<k,\\
q^{} {E}_{kj}{E}_{li}, & i=j<k<l, \\
q^{-1} {E}_{kj}{E}_{li}, & i<j<k=l,\\
{E}_{kj}{E}_{li} - \big(q-q^{-1}\big) {E}_{lj}{E}_{ki}, & i<j<l<k.
\end{cases}
\end{gather}
By the first and fourth lines above,
\begin{gather}\label{comm2}
\hat{E}_{il}\hat{E}_{jk} + q^{-1} \hat{E}_{ik}\hat{E}_{jl} = \hat{E}_{jk}\hat{E}_{il} + q \hat{E}_{jl}\hat{E}_{ik}, \\
\hat{E}_{li}\hat{E}_{kj} + q^{} \hat{E}_{ki}\hat{E}_{lj} = \hat{E}_{kj}\hat{E}_{li} + q^{-1} \hat{E}_{lj}\hat{E}_{ki}. \label{comm3}
\end{gather}

For any $ \ii $ and $ \jj $ in $ \mathbb{N}^m $, define
\[
\hat{E}^{\pm}_{\ii \jj} :=
\begin{cases}
\hat{E}_{i_1j_1}\cdots \hat{E}_{i_mj_m}, & \pm(i_1-j_1),\dots,\pm(i_m-j_m)\geq 0,\\
0, & \text{else}.
\end{cases}
\]
 and let $ e_{ \jj \ii }:= e_{j_1i_1} \otimes \dots \otimes e_{j_mi_m}$. For $\ii, \jj \in \mathcal{W}_m$, define
\[
\tilde{E}^{\pm}_{\ii\jj} =
\sum_{\zeta \in D_{\ii}}q^{\pm(\inv(\tau) -\inv(\zeta))}\hat{E}^{\pm}_{\bar{\zeta}(\ii)\bar{\tau}(\jj)},
\]
where $\tau$ is an arbitrary element of $D_{\jj}$. Note that the notation does not depend on $\tau$, which is justified by the following lemma and the relation $\inv(\bar{\tau}) = (m-1)m/2 - \inv(\tau)$.

\begin{Lemma} \label{prop 1}\quad
	\begin{enumerate}\itemsep=0pt
		\item[$1.$] For every $\tau,\tau' \in D_A$, the following identity holds: \begin{equation*}
		\sum_{\sigma\in D_B} q^{\mp (\inv(\sigma)-\inv(\tau))} \hat{E}^{\pm}_{\tau(A) \sigma (B)} = \sum_{\sigma\in D_B} q^{\mp (\inv(\sigma)-\inv(\tau'))} \hat{E}^{\pm}_{\tau'(A) \sigma (B)} .
		\end{equation*} 	 	
		\item[$2.$] Likewise, for every $\sigma,\sigma' \in D_B$, \begin{equation*}
		\sum_{\tau\in D_A} q^{\mp (\inv(\tau)-\inv(\sigma))} \hat{E}^{\pm}_{\tau(A) \sigma (B)} = \sum_{\tau\in D_A} q^{\mp (\inv(\tau)-\inv(\sigma'))} \hat{E}^{\pm}_{\tau(A) \sigma' (B)} .
		\end{equation*} 	 	
	\end{enumerate}
	\end{Lemma}
\begin{proof}
	For $m=2$, both cases are equivalent to the relations in \eqref{comm} and \eqref{comm1}.
	
	Now suppose that $ m>2$. We only prove part 1, as part 2 is similar. By Lemma \ref{Gene}, it suffices to consider the case when $\tau'=s\tau$ where $s$ is a transposition, and assume without loss of generality that $\inv(\tau')= \inv(\tau)+1$. Define the two sets $D_B^{(1)}$ and $D_B^{(2)}$ by
	\begin{gather*}
		D_B^{(1)}= \{ \sigma \in D_B \colon s \sigma(B)=B \},\qquad
		D_B^{(2)}= \{ \sigma \in D_B \colon s \sigma(B) \neq B \}.
	\end{gather*}	
	First, note that by the second and third lines of \eqref{comm},
\[
	\sum_{\sigma\in D_B^{(1)}} q^{\inv(\sigma)-\inv(\tau)} \hat{E}^-_{\tau(A) \sigma (B)} = \sum_{\sigma\in D_B^{(1)}} q^{\inv(\sigma)-\inv(s\tau)} \hat{E}^-_{s \tau(A) \sigma (B)} .
\]
	Partition the set $D_B^{(2)}$ into two sets $D_-$ and $D_+$ of equal cardinality, where $\sigma \in D_-$ if and only if $\inv(\sigma) < \inv(s \sigma)$. Then
	\begin{gather*}
	\sum_{\sigma\in D_B^{(2)}} q^{\inv(\sigma)-\inv(\tau)} \hat{E}^-_{\tau(A) \sigma (B)} = \sum_{\sigma \in D_-} \big( q^{\inv(\sigma)-\inv(\tau)} \hat{E}^-_{\tau(A) \sigma (B)} + q^{\inv(\sigma)+1-\inv(\tau)} \hat{E}^-_{\tau(A) s\sigma (B)} \big) \\
\qquad \stackrel{\eqref{comm2}}{=} \sum_{\sigma \in D_-} \big( q^{\inv(\sigma)+1-\inv(\tau')} \hat{E}^-_{s\tau(A) s\sigma (B)} + q^{\inv(\sigma)-1-\inv(\tau)} \hat{E}^-_{s\tau(A) \sigma (B)} \big)\\
\qquad =\sum_{\sigma \in D_-} \big( q^{\inv(\sigma)+1-\inv(\tau')} \hat{E}^-_{\tau'(A) s\sigma (B)} + q^{\inv(\sigma)-\inv(\tau')} \hat{E}^-_{\tau'(A) \sigma (B)}\big)\\
\qquad =\sum_{\sigma\in D_B^{(2)}} q^{\inv(\sigma)-\inv(\tau')} \hat{E}^-_{\tau'(A) \sigma (B)} ,
	\end{gather*}
	as needed.
	The proof for $E^+$ is similar, where one uses \eqref{comm3} instead of \eqref{comm2}.
\end{proof}

\begin{Example}
Consider $m=2$ and $N=3$. Set $\jj=(0,1)$ and $\ii=(2,3)$. Then
\[
\tilde{E}_{\jj\ii}^- = \hat{E}_{13}\hat{E}_{02} + q \hat{E}_{03}\hat{E}_{12} = q^{-1}\big(\hat{E}_{12}\hat{E}_{03} + q \hat{E}_{02}\hat{E}_{13} \big),
\]
with the equality following from \eqref{comm2}. Additionally,
\[
\tilde{E}_{\ii\jj}^+ = \hat{E}_{31}\hat{E}_{20} + q^{-1}\hat{E}_{21}\hat{E}_{30} = q\big( \hat{E}_{30}\hat{E}_{21} + q^{-1}\hat{E}_{20}\hat{E}_{31}\big),
\]
with the equality following from \eqref{comm3}.
\end{Example}

\section{Statements and proofs}

The main theorem is the following expression for central elements of $\Uqgl$.
\begin{Theorem}\label{Main}
	The element given by
\begin{equation*}
C_m=\sum_{\ii,\jj \in \mathcal{W}_m} q^{ 2i_1 + \dots + 2i_m-Nm} \tilde{E}_{\jj \ii}^-\tilde{E}_{\ii\jj}^+
\end{equation*}
is central in $ \Uqgl $.
\end{Theorem}

\begin{Example}
Consider $m=1$. Then
\[
C_1 = \big(q-q^{-1}\big)^{-2}\sum_{i=0}^N q^{2i-N} q^{2\epsilon_i} + \sum_{0\leq j<i \leq N} q^{2i-N-1} q^{2\epsilon_i+2\epsilon_j}E_{ji}E_{ij},
\]
which is (up to a constant) the central element $C$ from \cite{GZB}.
\end{Example}
\begin{Example}
Consider $m=2$ and $N=3$. Then
\begin{gather*}
C_2 = \sum_{0 \leq j \leq i \leq 3} q^{-6}q^{4i}\hat{E}_{ji}^2 \hat{E}_{ij}^2 \\
\hphantom{C_2 =}{}
+ \sum_{0\leq j \leq i_1 < i_2 \leq 3} q^{-6+2i_1+2i_2}\hat{E}_{ji_1}\hat{E}_{ji_2}\big(\hat{E}_{i_2j}\hat{E}_{i_1j} + q^{-1}\hat{E}_{i_1j}\hat{E}_{i_2j}\big) \\
\hphantom{C_2 =}{}+ \sum_{0\leq j_1 < j_2 \leq i \leq 3} q^{-6+4i}\big(\hat{E}_{j_1i}\hat{E}_{j_2i} + q^{-1}\hat{E}_{j_2i}\hat{E}_{j_1i} \big) \hat{E}_{ij_2}\hat{E}_{ij_1}\\
\hphantom{C_2 =}{}+ \sum_{0\leq j_1\leq i_1 < j_2 \leq i_2 \leq 3} q^{-6+2i_1+2i_2} \hat{E}_{j_1i_1}\hat{E}_{j_2i_2}\hat{E}_{i_2j_2}\hat{E}_{i_1j_1}\\
\hphantom{C_2 =}{}+\sum_{0\leq j_1 < j_2\leq i_1 < i_2 \leq 3}\!\! q^{-6+2i_1+2i_2} \big(\hat{E}_{j_1i_1}\hat{E}_{j_2i_2} + q^{-1}\hat{E}_{j_ii_2}\hat{E}_{j_2i_1} \big)\big(\hat{E}_{i_2j_2}\hat{E}_{i_1j_1} + q^{-1}\hat{E}_{i_2j_1}\hat{E}_{i_1j_2}\big).
\end{gather*}
The central element $C_2$ has $50$ terms, consisting of $10$ terms from $0\leq j \leq i \leq 3$, $20$ terms from $0 \leq j_1<j_2 \leq i \leq 3$ and $0\leq j\leq i_1<i_2 \leq 3$, $15$ terms coming from $0\leq j_1\leq i_1 < j_2 \leq i_2 \leq 3$ and $5$ terms coming from $0\leq j_1 < j_2\leq i_1 < i_2 \leq 3$. One can also verify that~$50$ is the correct number of terms, from the fact that $|\mathcal{W}_2|=10$ and the set $\{ \{\ii,\jj\}\colon \ii,\jj \in \mathcal{W}_2\}$ has $55=10(10+1)/2$ elements, but the term $\tilde{E}^{\pm}_{\jj\ii}$ is zero when $\{\ii,\jj\}$ is one of the $5$ sets $\{(0,2),(1,1)\},\{(0,3),(2,2)\},\{(0,3),(1,1)\},\{(1,3),(2,2)\},\{(1,2),(0,3)\}$.

The representation $P_2 \mathbb{C}^4$ of $\Uqgl$ is $10$-dimensional, which can be explicitly written using Example \ref{ExRep} below. By multiplying $10\times 10$ matrices, one can check that $C_2$ acts as $\text{const}\cdot \mathrm{Id}_{10}$, where the constant is $\big(q+q^{-1}\big)^{-4}q^{-6}\big(1+q^2+2q^4+q^6+q^8+q^{10}+q^{12}+q^{14}+q^{20}\big)$.
\end{Example}

\begin{Remark}\label{Remark1}
In \cite{ZGBJPhysA}, the central element $C^{\Lambda_0}$ is defined, where $\Lambda_0$ is the highest weight of a finite-dimensional irreducible module $V(\Lambda_0)$. The construction there is similar to the one here, with the major difference being the use of explicit universal $R$-matrices in place of fused $R$-matrices. Although it is not necessarily simple to check directly that $C_m$ equals (up to a~constant) $C^{\Lambda_0}$ for $V(\Lambda_0)=P_mV^{\otimes m}$, it is straightforward to check that their eigenvalues are the same.

If $V(\Lambda_0)$ has distinct weights $\lambda_1,\dots,\lambda_r$ with multiplicities $d_1,\dots,d_r$, then the eigenvalue of~$C^{\Lambda_0}$ on an irreducible module with highest weight $\Lambda$ is given by
\[
\sum_{k=1}^r d_k q^{2(\Lambda + \rho,\lambda_k)}.
\]
Here, $\rho$ is half the sum of the positive roots, and $(\cdot,\cdot)$ is the usual invariant bilinear form on $\mathfrak{h}^*$. Now take $\Lambda_0$ to be the highest weight of $P_mV^{\otimes m}$; then the distinct weights are elements of $\mathcal{B}^{(N)}_m$ with multiplicity $1$. Therefore the eigenvalue is
\[
\sum_{\mu \in \mathcal{B}^{(N)}_m} q^{(2\rho,\mu)}q^{2(\Lambda,\mu)}.
\]

If $C_m$ acts on the same module $V(\Lambda)$, then its eigenvalue can be found by evaluating on the lowest weight vector, because then only the diagonal terms $(\ii=\jj)$ have a nonzero contribution. So the eigenvalue is
\[
\big(q-q^{-1}\big)^{-2m}q^{-Nm}\sum_{\ii \in \mathcal{W}_m} q^{ 2i_1 + \dots + 2i_m} q^{(2\mu^{(N)}(\ii),\Lambda)}.
\]
By Lemma \ref{Leem}, this is $\big(q-q^{-1}\big)^{-2m}q^{-Nm}$ times the eigenvalue of $C^{\Lambda_0}$.
\end{Remark}

\begin{Remark}
 The value of the central character on $C_m$ can be found using
\cite[Corollary~A.2]{mudrov2007quantum}. Their result implies that
\[
\chi^{\lambda}(C_m) = \operatorname{Tr}\big(\pi_{W_m}\big(q^{2(h_{\lambda}+h_{\rho})}\big)\big),
\]
where $W_m$ is the $m$th symmetric tensor power of the defining module. Since $W_m$ has weights indexed by $\mathcal{B}_m^{(N)}$, the right-hand side is
\[
\sum_{\nu \in \mathcal{B}_m^{(N)}} \prod_{i=0}^Nq^{ (2\lambda_i + N-2i)\nu_i }.
\]
\end{Remark}

\subsection[Basis for P\_mV\^\{otimes m\}]{Basis for $\boldsymbol{P_mV^{\otimes m}}$}

Before proving Theorem \ref{Main}, we will write a basis for the symmetric projection $P_mV^{\otimes m}$.

Let $ I_0,\dots,I_N $ be the canonical basis of $ V $, and define the action on $ V $ by 	
\begin{gather}
\hat{e}_{+,i}I_j=\delta_{ij}I_{j-1},\qquad
\hat{e}_{-,i}I_{j-1}=\delta_{ij}I_j,\qquad
q^{\pm \epsilon_i/2}I_j=q^{\pm \delta_ij/2}I_j.\label{21}
\end{gather}
For $\mu \in \mathcal{B}_m^{(N)}$, define the vector
\begin{equation*}
v\zo =I_0^{\otimes \mu_0} \otimes \dots \otimes I_N^{\otimes \mu_N}
\end{equation*}
Define
\begin{equation*}
M(\zz) = \sum_{\sigma \in D^{\mu}}q^{-\inv (\sigma)} \sigma (v(\mu)).
\end{equation*}
Here, as before, $S_m $ acts on $ V^{\otimes m} $ by permuting the order.
By an abuse of notation, for $\ii \in \mathcal{W}_m$ we define
\[
v(\ii) = v\big(\mu^{(N)}(\ii)\big), \qquad M(\ii) = M\big(\mu^{(N)}(\ii)\big),
\]
where $N \geq i_m$. Note that these definitions do not depend on the value of $N$.
We briefly note the identity
\begin{equation}\label{Brief}
e^{\sigma(a)}_{\pm,i} ( \sigma(v(\mu)) ) = \sigma\big( e^{(a)}_{\pm,i} ( v(\mu) ) \big) \qquad \text{for all } \sigma \in S_m.
\end{equation}

We now show that the set $ \{M(\mu)\}_{\mu\in \mathcal{B}^{(N)}_m} $ gives a basis for $ P_{m}V^{\otimes m} $. A more general statement appeared in~\cite{KIMRN} with a more complicated proof, but the expression here is more convenient for calculations.
\begin{Theorem}
	The set $ \{M(\mu)\}_{\mu\in \mathcal{B}^{(N)}_m} $ gives a basis for $ P_{m}V^{\otimes m} $.
\end{Theorem}
\begin{proof}
Note that $ \big|\mathcal{B}^{(N)}_m\big|= \dimension P_mV^{\otimes m} $ and $ \{M(\mu)\}_{\mu\in \mathcal{B}^{(N)}_m} $ is a linearly independent set, so it suffices to show that $ M(\mu) \in P_mV^{\otimes m} $ for all $ \mu \in \mathcal{B}^{(N)}_m $. To show that every $ M(\mu) $ is in $ P_{m}V^{\otimes m} $, it suffices to show that \begin{equation} \label{hati}
\Delta^{(m)}(\hat{e}_{-,i}) (M(\mu))=q^{\frac{1}{2}(\mu_{i-1}+\mu_i-1)} \big(1+q^{-2}+q^{-4}+\dots+q^{-2\mu_i}\big) \cdot M\big(\mu-\hat{i}\big),
\end{equation}
where $ \hat{i}=(0,\dots,0,1,-1,0,\dots,0) $ with the $ -1 $ occurring in the $ i $th position. By~\eqref{qerel} and~\eqref{21}, $ \Delta^{(m)} (\hat{e}_{-,i})(M(\mu)) $ is in the span of $ \{\tau(v(\mu-\hat{i}))\}_{\tau\in S_m}$. Let $ A(\tau) $ be the coefficients in the expansion
\[
 \Delta^{(m)} (\hat{e}_{-,i})(M(\mu)) = \sum_{\tau \in D^{\mu-\hat{i}}} A(\tau)\tau \big(v\big(\mu-\hat{i}\big)\big).
\]
It suffices to show that $A(\tau)= q^{\frac{1}{2}(\mu_{i-1}+\mu_i-1)}\big(1+q^{-2}+q^{-4}+\dots+q^{-2\mu_i}\big) q^{-\inv(\tau)}$ for all $\tau \in D^{\mu-\hat{i}}$. In fact, we show something stronger: for all $\tau \in D^{\mu-\hat{i}}$, there exist elements $\sigma^{(0)},\dots, \sigma^{(\mu_i)}\in D^{\mu}$ such that for $0 \leq j \leq \mu_i$,
\[
\hat{e}_{-,i}^{(a_j)} \big( q^{-\inv(\sigma^{(j)})} \sigma^{(j)}(v(\mu)) \big) = q^{\frac{1}{2}(\mu_{i-1}+\mu_i-1)} q^{-2j} q^{-\inv(\tau)} \tau\big(v\big(\mu- \hat{i}\big)\big).
\]

We will proceed by induction on the value of $\inv(\tau)$, using Lemma \ref{Gene}. The base case is when~$\tau$ is the identity permutation. Then it is straightforward to check that
\[
\sigma^{(j)} = (\mu_{[0,i-1]} \ \ \mu_{[0,i-1]}+1 \ \ \dots \ \ \mu_{[0,i-1]}+j),
\]
where $\mu_{[0,i-1]}= \mu_0 + \dots + \mu_{i-1}$, satisfies the necessary conditions.

Now fix $\tau \in D^{\mu - \hat{i}}$, and suppose that the induction hypothesis holds for $\tau$. Suppose that $\tilde{\tau} \in D^{\mu - \hat{i}}$ satisfies $\tilde{\tau}=s\tau$ for some transposition $s$, and assume without loss of generality that $\inv(\tilde{\tau}) = \inv(\tau) + 1$. Define
\[
\tilde{\sigma}^{(j)}
=
\begin{cases}
s \sigma^{(j)}, & \text{if } s \sigma^{(j)} \in D^{\mu}, \\
\sigma^{(j)}, & \text{else. }
\end{cases}
\]
We now aim to prove that
\begin{equation}\label{Aim}
\hat{e}^{(s(a_j))}_{-,i}\big(q^{-\inv(\tilde{\sigma}^{(j)})} \tilde{\sigma}^{(j)}(v(\mu)) \big) = q^{\frac{1}{2}(\mu_{i-1}+\mu_i-1)} q^{-2j} q^{-\inv(\tilde{\tau})} \tilde{\tau}\big(v\big(\mu- \hat{i}\big)\big).
\end{equation}
If $\tilde{\sigma}^{(j)} = s\sigma^{(j)}$, then \eqref{Aim} follows from \eqref{Brief} and the induction hypothesis. Now assume that $\tilde{\sigma}^{(j)} = \sigma^{(j)}$. Then $s= (a_j-1 \ \ a_j)$ and
\[
\tilde{\sigma}^{(j)} ( v(\mu) ) = \underbrace{I_* \otimes \dots \otimes I_* }_{a_j-2} \otimes I_{i-1} \otimes I_{i-1} \otimes \underbrace{I_* \otimes \dots \otimes I_*}_{m-a_j}.
\]
Using that
\[
 P \big( \big(\hat{e}_{i,-} \otimes q^{-h_i/2}\big)(I_{i-1} \otimes I_{i-1}) \big) = q^{-1}\big(q^{h_i/2} \otimes \hat{e}_{i,-}\big)(I_{i-1} \otimes I_{i-1}) ,
\]
where $P(x\otimes y)=y\otimes x$ is the permutation operator, we have that
\[
\hat{e}^{(s(a_j))}_{-,i}\big(q^{-\inv(\tilde{\sigma}^{(j)})} \tilde{\sigma}^{(j)}(v(\mu)) \big) = q^{-1}\hat{e}^{(a_j)}_{-,i}\big(q^{-\inv({\sigma}^{(j)})} {\sigma}^{(j)}(v(\mu)) \big).
\]
And now \eqref{Aim} follows from the induction hypothesis.
\end{proof}
\begin{Remark}\label{ExRep}
	For all $ \mu \in \mathcal{B}^{(N)}_m $, let $ \tilde{M}(\mu)=c(\mu)M(\mu) $ where $ c(\mu) $ is defined inductively by
	\begin{equation*}
	c(m,0,\dots,0)=1, \qquad
	\frac{c(\mu-i)}{c(\mu)}=\frac{q^{\frac{\mu_{i-1}+\mu_i+1}{2}}\big(1-q^{-2\mu_i}\big)}{q^{\mu_i+1}-q^{-\mu_i-1}} ,
	\end{equation*} then equation \eqref{hati} is equivalent to
\begin{equation*}
	\Delta^{(m)} (\hat{e}_{-,i}) \tilde{M}(\mu) =\frac{q^{\mu_i+1}-q^{-\mu_i-1}}{q-q^{-1}} \tilde{M}\big(\mu-\hat{i}\big).
	\end{equation*} In fact, one can show that
	\begin{gather*}
		\Delta^{(m)} (\hat{e}_{+,i}) \tilde{M}(\mu) =\frac{q^{\mu_{i-1}+1}-q^{-\mu_{i-1}-1}}{q-q^{-1}} \tilde{M}\big(\mu+\hat{i}\big),\\
		\Delta^{(m)} (\hat{e}_{-,i}) \tilde{M}(\mu) =\frac{q^{\mu_i+1}-q^{-\mu_i-1}}{q-q^{-1}} \tilde{M}\big(\mu-\hat{i}\big),\\
	 \Delta^{(m)} \big(q^{\pm h_i/2}\big) \tilde{M}(\mu) =q^{\pm\frac{1}{2}(\mu_{i-1}-\mu_i)} \tilde{M}(\mu)
	\end{gather*} defines a representation on $ P_mV^{\otimes m} $. This is equivalent to the representation in \cite[equation~(3)]{KMMO} and \cite[Lemma~3.1]{KIMRN}.
\end{Remark}
\begin{Corollary}\label{Core} Let $\ii,\jj\in\mathcal{W}_m$. For any $\tau,\sigma \in D_{\ii}$ and any $\zeta \in D_{\jj}$,
\[
\big(q^{\inv(\tau)}e_{\zeta(\jj) \tau(\ii)} - q^{\inv(\sigma)}e_{\zeta(\jj) \sigma(\ii)} \big) \big|_{P_{ m}^+V^{\otimes m} } =0.
\]
	
	Furthermore,
\[
\bigg( \sum_{\tau \in D_{\jj}} q^{-d_{\ii}(\tau)}e_{\tau(\jj)\tau(\ii)} \bigg) M(\ii) = M(\jj) .
\]
\end{Corollary}
\begin{proof}
For any $\mu \in \mathcal{B}^{(N)}_m$ not equal to $\mu(\ii)$, it is straightforward that
\[
\big(q^{\inv(\tau)}e_{\zeta(\jj) \tau(\ii)} - q^{\inv(\sigma)}e_{\zeta(\jj) \sigma(\ii)} \big) M(\mu)=0-0=0.
\]
On the other hand,
\begin{gather*}
\begin{split}
&\big(q^{\inv(\tau)}e_{\zeta(\jj) \tau(\ii)} - q^{\inv(\sigma)}e_{\zeta(\jj) \sigma(\ii)} \big) M(\mu(\ii))\\
&\qquad = q^{\inv(\tau)}e_{\zeta(\jj) \tau(\ii)} \big(q^{-\inv(\tau)} \tau(v(\mu(\ii)) )\big) - q^{\inv(\sigma)}e_{\zeta(\jj) \sigma(\ii)} \big(q^{-\inv(\sigma)} \sigma(v(\sigma(\ii)) )\big) \\
&\qquad = \zeta(v(\mu(\jj))) - \zeta(v(\mu(\jj)))=0.
\end{split}
\end{gather*}

For the second statement, we use that
\[
e_{\tau(\jj)\tau(\ii)}\sigma(v(\ii))
=
\begin{cases}
\tau(v(\jj)), & \tau \sim_{\ii} \sigma, \\
0, & \text{else.}
\end{cases}
\]
in order to show
\begin{align*}
\bigg( \sum_{\tau \in D_{\jj}} q^{-d_{\ii}(\tau)}e_{\tau(\jj)\tau(\ii)} \bigg) \sum_{\sigma \in D_{\ii}} q^{-\inv(\sigma)} \sigma(v(\ii)) &= \sum_{\substack{\tau \in D_{\jj},\sigma \in D_{\ii} \\ \tau\sim_{\ii} \sigma}}q^{-d_{\ii}(\tau)+ \inv(\sigma)} e_{\tau(\jj)\tau(\ii)}\sigma(v(\ii))\\
&= \sum_{\tau \in D_{\jj}} q^{-\inv(\tau)} \tau(v(\jj)).\tag*{\qed}
\end{align*}\renewcommand{\qed}{}
\end{proof}

\begin{Example}
Some examples are
\begin{gather*}
\big(e_{21} \otimes e_{31} + q^{-1}e_{31} \otimes e_{21}\big)M(\mathbf{1,1})=M(\mathbf{2,3}),\\
(e_{31} \otimes e_{32})M(\mathbf{1,2}) = M(\mathbf{3,3}),\\
(e_{31} \otimes e_{42} + e_{42} \otimes e_{31} )M(\mathbf{1,2}) = M(\mathbf{3,4}).
\end{gather*}
\end{Example}

With Corollary \ref{Core} as motivation, define for any $\ii,\jj \in \mathcal{W}_m$,
\[
\tilde{e}_{\jj \ii} = \sum_{\tau \in D_{\jj}} q^{-d_{\ii}(\tau)}e_{\tau(\jj)\tau(\ii)} .
\]

\subsection{Proof of Theorem \ref{Main}}
We will now prove Theorem \ref{Main}. Begin by rewriting $R_{0,\{1,\dots,m\}}$. By definition,
\begin{align*}
R_{0,\{1,\dots,m\}} &= R_{0m}R_{0m-1}\cdots R_{01} \\
&=\bigg( \sum_{i_m \geq j_m} \hat{E}_{i_mj_m} \otimes \id^{\otimes m-1} \otimes e_{j_mi_m} \bigg) \cdots \bigg(\sum_{i_1\geq j_1}\hat{E}_{i_1j_1}\otimes e_{j_1i_1} \otimes \id^{\otimes m-1}\bigg)
\\&= \sum_{i_m \geq j_m} \cdots \sum_{i_1\geq j_1} \hat{E}_{i_mj_m} \cdots \hat{E}_{i_1j_1} \otimes e_{j_1i_1} \otimes \dots \otimes e_{j_mi_m} \\
&= \sum_{\ii,\jj \in \mathcal{W}_m } \sum_{\zeta \in D_{\ii},\tau \in D_{\jj}}\hat{E}_{\bar{\zeta}(\ii)\bar{\tau}(\jj)}^{+} \otimes e_{\tau(\jj)\zeta(\ii)}.
\end{align*}
From here, the goal is to write this expression as an element of $\Uqgl \otimes \End\big(P_m\big(V^{\otimes m}\big)\big)$. By Corollary \ref{Core},
\[
R_{0,\{1,\dots,m\}} = \sum_{\ii,\jj \in \mathcal{W}_m } \sum_{\zeta \in D_{\ii},\tau \in D_{\jj}}q^{\inv(\sigma_{\ii}(\tau)) -\inv(\zeta)}\hat{E}^+_{\bar{\zeta}(\ii)\bar{\tau}(\jj)} \otimes e_{\tau(\jj)\sigma_{\ii}(\tau)(\ii)}.
\]
By definition, $\tau = \sigma_{\ii}(\tau) \xi_{\ii}(\tau)$ and $d_{\ii}(\tau) = \inv(\xi_{\ii}(\tau))$, so therefore
\[
R_{0,\{1,\dots,m\}} = \sum_{\ii,\jj \in \mathcal{W}_m } \sum_{\zeta \in D_{\ii},\tau \in D_{\jj}}q^{\inv(\tau) -d_{\ii}(\tau) -\inv(\zeta)}\hat{E}^+_{\bar{\zeta}(\ii)\bar{\tau}(\jj)} \otimes e_{\tau(\jj)\tau(\ii)}.
\]
By Lemma \ref{prop 1} and the identity $\inv(\bar{\tau}) = m(m-1)/2 - \inv(\tau)$,
\[
\sum_{\zeta \in D_{\ii}}q^{\inv(\tau) -\inv(\zeta)}\hat{E}_{\bar{\zeta}(\ii)\bar{\tau}(\jj)}
\]
does not depend on $\tau$. Therefore, the $R$-matrix factors as
\begin{align*}
R_{0,\{1,\dots,m\} } &= \sum_{\ii,\jj \in \mathcal{W}_m } \bigg( \sum_{\zeta \in D_{\ii}}q^{\inv(\tau) -\inv(\zeta)}\hat{E}^+_{\bar{\zeta}(\ii)\bar{\tau}(\jj)} \bigg) \otimes \bigg( \sum_{\tau \in D_{\jj}} q^{-d_{\ii}(\tau)}e_{\tau(\jj)\tau(\ii)} \bigg) \\
&= \sum_{\ii,\jj \in \mathcal{W}_m } \tilde{E}_{\ii\jj}^+ \otimes e_{\jj\ii}.
\end{align*}
A similar argument shows that
\[
R^T_{0,\{1,\dots,m\}}=\sum_{\ii,\jj \in \mathcal{W}_m } \tilde{E}^-_{\jj\ii} \otimes e_{\ii\jj}.
\]

Therefore,
\[
\Gamma_m = R^T_{0,\{1,\dots,m\} } R_{0,\{1,\dots,m\} } = \sum_{\ii,\jj \in \mathcal{W}_m } \sum_{\ii',\jj' \in \mathcal{W}_m } \tilde{E}_{\jj \ii}^-\tilde{E}_{\ii'\jj'}^+ \otimes e_{\ii\jj}e_{\jj'\ii'}.
\]
By Corollary \ref{Core},
\[
C_m = (\id \otimes \tr_q)(\Gamma_m) = \sum_{\ii,\jj \in \mathcal{W}_m } \tilde{E}_{\jj \ii}^-\tilde{E}_{\ii\jj}^+ \otimes \tr_q(e_{\ii\jj}e_{\jj\ii}).
\]

It just remains to calculate the quantum trace. It is given by
\[
\tr_q(e_{\ii\jj}e_{\jj\ii}) = \tr_q(e_{\ii\ii}) = \tr\big(q^{-NE_{00} -(N-2)E_{11} + \dots + NE_{NN}} e_{\ii\ii}\big) = q^{-N\mu_{0} -(N-2)\mu_{1} + \dots + N\mu_{N}},
\]
where $\mu=\mu^{(N)}(\ii)$. Applying Lemma \ref{Leem} finishes the proof.

\subsection*{Acknowledgements}
Jeffrey Kuan was supported by the Minerva Foundation and NSF grant DMS-1502665. Both authors were supported by the 2017 Columbia Mathematics REU program, which was funded by Columbia University. We thank Yi Sun for helpful discussions. We thank an anonymous referee for providing references for the introduction.

\pdfbookmark[1]{References}{ref}
\LastPageEnding

\end{document}